\documentclass[a4paper,reqno]{amsart}
\usepackage{latexsym}
\usepackage{amssymb,amsthm,amsmath}
\usepackage{graphicx}

\usepackage{geometry}
\usepackage{eurosym}%
\usepackage{pgfgantt}
\usepackage{pdfpages}
\usepackage[pdfborder={0 0 0}]{hyperref}
\usepackage{bigints}

\usepackage{setspace}
\calclayout

\theoremstyle{plain}
\newtheorem{theorem}{Theorem}
\newtheorem{lemma}[theorem]{Lemma}
\newtheorem{corollary}[theorem]{Corollary}

\theoremstyle{definition}

\theoremstyle{remark}

\def\d#1{{#1\kern-0.4em\char"16\kern-0.1em}}
\def\D#1{{\raise0.2ex\hbox{-}\kern-0.4em #1}}
\reversemarginpar

\newcounter{zd}

\newcounter{zdr}[subsection]

\def\mx{{\bf x}}

\def\A{{\cal A}}

\def\cal{\mathcal}
\let\mib=\boldsymbol

\def\R{{\bf R}}

\def\N{{\bf N}}

\def\mxi{{\mib \xi}}

\def\mx{{\bf x}}

\def\PP{\mathbb{P}}

\def\u{{\bf u}}

\def\w{{\bf w}}

\begin{document}

\title[High-frequency tail diagnostics]{A high-frequency tail condition and a diagnostic iteration for the Navier--Stokes equations}

\author{D.~Mitrovi\'{c}}\address{Darko Mitrovi\'{c}, Faculty of Mathematics,
University of Vienna, Oskar Morgenstern-Platz 1,
1090 Vienna, Austria and } \address{Faculty of Mathematics and Natural Sciences, University of Montenegro, George Washington boulevard BB, 81000 Podgorica, Montenegro} \email{darkom@ucg.ac.me}

\subjclass[2020]{35Q30, 76D05, 76D03, 35B44}

\keywords{Navier--Stokes equations; Leray solutions; frequency localization; high--frequency cutoff; diagnostic Picard iteration; blow-up exclusion}

\begin{abstract}
We consider Leray solutions of the three--dimensional incompressible Navier--Stokes
equations on $\R^3$ with smooth, rapidly decaying initial data.  The analysis is
based on a frequency decomposition into low and high modes via the cutoffs
$\A_R=\phi(|D|/R)$ and $\A^R=I-\A_R$.

Combining the energy inequality with Bernstein estimates yields uniform control of
the low--frequency component $\A_R\u$.  For the high--frequency component we assume a
quantitative \emph{turbulence condition}, requiring that the solution possesses a
non--negligible high--frequency tail in $L^\infty$ (in fact, it suffices to impose
this condition only on a terminal time layer near a putative blow--up time).

Under this hypothesis we introduce a time--localized diagnostic Picard iteration
adapted to $\A^R\u$.  Using a uniform $L^\infty$ estimate of Giga--Inui--Matsui type
(with the cutoff $\A^R$) together with high--frequency heat--flow decay, we show that
the iteration is contractive and converges to $\A^R\u$, providing a uniform bound for
$\A^R\u$ up to the maximal time of boundedness.  Consequently, the turbulence regime
is incompatible with finite--time blow--up: any Leray solution satisfying the
turbulence condition is bounded, and hence smooth, for all times (equivalently, it cannot blow up in finite time).
\end{abstract}

\maketitle

\section{Introduction}\label{sec:intro}

The three--dimensional incompressible Navier--Stokes equations are among the most
intensively studied systems in mathematics and mathematical physics.  Despite more
than a century of work---and despite the existence of a robust global theory of
finite--energy weak solutions---the fundamental question of whether smooth solutions
can develop finite--time singularities remains open.  The Millennium Prize Problem
formulation highlights precisely this tension: the equations are simple to state, yet
their global regularity theory in dimension three is still incomplete.

The modern mathematical theory encompasses a wide range of phenomena, from the
classical construction of global weak solutions due to Leray to striking recent
developments showing non--uniqueness and instability at low regularity.  We do not
attempt to survey this vast literature and instead refer to \cite{rieusset} for a
comprehensive overview.  Among representative advances, let us mention results on
non--uniqueness of weak or Leray--type solutions \cite{colombo,vicol,arxiv1,arxiv2},
which underscore the subtlety of the solution concept, as well as contributions to
regularity and conditional uniqueness, including the classical Serrin criteria
\cite{serrin}, the partial regularity theory of Caffarelli--Kohn--Nirenberg \cite{caf},
approaches in critical or borderline spaces such as $L^{3,\infty}$ \cite{ISS}, and
refined local--in--time and initial--time regularity results \cite{giga,jia}.

\smallskip

Let us recall that the incompressible Navier--Stokes equations on $\R^3$ read
\begin{equation}\label{N-S}
\begin{aligned}
&\partial_t \u - \Delta \u + (\u\cdot\nabla)\u + \nabla p = 0,\\
&\nabla\cdot \u = 0,\\
&\u|_{t=0} = \u_0,\qquad \nabla\cdot \u_0 = 0,
\end{aligned}
\end{equation}
where $\u_0$ is smooth and rapidly decaying (for definiteness, one may take
$\u_0\in\mathcal{S}(\R^3)$).  By incompressibility,
\[
(\u\cdot\nabla)\u = \nabla\!\cdot(\u\otimes\u).
\]
The pressure admits the representation
\begin{equation}\label{pressure}
p = \sum_{k,j=1}^3 R_kR_j(u_ku_j),
\end{equation}
where $R_j$ denote the Riesz transforms, i.e.\ Fourier multipliers with symbols
$i\,\xi_j/|\mxi|$.  Consequently, \eqref{N-S} can be written in the projected form
\begin{equation}\label{N-S-L}
\partial_t \u + \mathbb{P}\,\nabla\!\cdot(\u\otimes\u)=\Delta\u,
\end{equation}
where $\mathbb{P}$ is the Leray projector onto divergence--free vector fields.

\smallskip

{\em Motivation: high frequencies and ``turbulence'' as a power--law tail regime.}
A recurring heuristic in the discussion of possible singularity formation is that
nonlinear effects might drive a sustained transfer of activity toward smaller and
smaller spatial scales (equivalently, toward higher and higher frequencies), until the
linear dissipation can no longer compensate.  From a purely mathematical point of
view, an $L^\infty$--blow--up scenario must indeed involve high frequencies: if for
some fixed $R$ one had $\A^R\u(t)\equiv 0$ on a time interval, then $\u(t)$ would be
band--limited and Bernstein inequalities would prevent $\|\u(t)\|_\infty$ from becoming
unbounded.  Thus, any genuine singular behavior must be accompanied by nontrivial
high--frequency content for every fixed cutoff.

In this paper we encode a quantitative \emph{power--law non--negligibility} of the
high--frequency tail.  Fix $R\gg 1$ and decompose the solution into low and high modes
via
\begin{equation}\label{AR}
\A_R := \phi(|D|/R),\qquad \A^R:=I-\A_R,
\end{equation}
where $\phi\in C_c^\infty(\R)$ is nonnegative, $0\le \phi\le 1$, supported in $(-2,2)$,
and equal to $1$ on $(-1,1)$.  Accordingly,
\[
\u=\A_R\u+\A^R\u .
\]
We say that the solution is \emph{turbulent (at scale $R$)} on a time interval if,
for some parameters $k\geq 2$, $p\in[1,\infty]$ and $\epsilon>0$, and for $R$ large
enough, one has
\begin{equation}\label{turbulent}
\|\u(t)\|_\infty
\le
\frac{R^{\frac{k-3/p-\epsilon}{2}}}{\|\nabla^k \u_0\|^{1/2}_{p}}\,
\|\A^R\u(t)\|_\infty.
\end{equation}
Here and throughout,
\[
\|\A^R\u(t)\|_\infty := \|\A^R\u(t,\cdot)\|_{L^\infty(\R^3)}.
\]
Condition \eqref{turbulent} should not be read as a literal domination of $\A^R\u$
over $\u$: the prefactor $R^{\frac{k-3/p-\epsilon}{2}}$ allows the ratio
$\|\u(t)\|_\infty/\|\A^R\u(t)\|_\infty$ to grow with $R$.
Rather, \eqref{turbulent} enforces that the high--frequency tail cannot be
\emph{too small} relative to the full amplitude at scale $R$---it is permitted to be
smaller, but only down to a prescribed power--law threshold.  In particular, the
high--frequency component is ruled out from being super--polynomially negligible as
$R\to\infty$, which is consistent with the intuition that ``wild'' behavior requires
persistent activity at small scales.  This viewpoint is also compatible with the
``ascending chain'' philosophy used in the study of turbulent scenarios, where
power--law relations between higher and lower derivatives are interpreted as evidence
of significant high--frequency content (see e.g. \cite{GX}).

{\em Main result.}
Our main theorem shows that, within the Leray framework, the turbulent tail condition \eqref{turbulent} prevents finite--time blow--up.

\begin{theorem}\label{main}
Any Leray solution to the Navier--Stokes equations \eqref{N-S} corresponding to
divergence--free Schwartz initial data $\u_0$ and satisfying \eqref{turbulent} for $R$
large enough is bounded (and hence does not blow up).
\end{theorem}

In Section 3 we refine this and show it suffices to assume \eqref{turbulent} only on a terminal time layer near a putative blow-up time.

\smallskip

{\em Idea of the proof.}
The argument is based on rewriting the mild formulation in terms of the high--frequency
component $\A^R\u$.  We treat $\u$ as a given function satisfying \eqref{turbulent} and
construct, via a Picard iteration, a bounded sequence converging to $\A^R\u$ in
$L^\infty$.  The key point is that the effective nonlinearity in the $\A^R$--equation
is weakened by the cutoff and can be controlled by combining Bernstein--type bounds
with the local $L^\infty$ theory.

\smallskip

{\em Comparison with Tao's averaged Navier--Stokes model.}
A particularly relevant point of comparison is provided by the averaged Navier--Stokes
model introduced by Tao~\cite{Tao}.  In that work, the nonlinear term is replaced by an
averaged bilinear form of the schematic type
\[
\widetilde{B}(u,u)=m_3(D)\,\langle m_2(D)u,\; m_1(D)u\rangle,
\]
where $m_1,m_2,m_3$ are zero--order H\"ormander--Mikhlin multipliers.  Tao proves that
for this modified system there exist smooth rapidly decaying initial data whose
corresponding solutions have non--increasing energy and yet blow up in finite time,
while the averaged nonlinearity preserves the formal cancellation properties needed to
construct global Leray--type weak solutions.

From the perspective adopted here, Tao's construction highlights a precise obstruction
to $L^\infty$--based control of the nonlinear term.  Inserting $\widetilde{B}$ into the
mild formulation and formally taking the $L^\infty$ norm yields
\begin{equation}\label{giga-28-org-nogo}
\bigl\| \nabla\!\cdot e^{t\Delta} m_3(D)\,
      \langle m_2(D)u,\; m_1(D)u\rangle \bigr\|_{L^\infty}
\le C_G\,t^{-1/2}\,
      \|\langle m_2(D)u,\; m_1(D)u\rangle\|_{L^\infty},
\end{equation}
where $m_3(D)$ drops out in the same way as the Leray projector in the classical
estimate.  The difficulty lies in the right--hand side: the additional multipliers
$m_1(D)$ and $m_2(D)$ prevent direct $L^\infty$ control of
$\langle m_2(D)u,\;m_1(D)u\rangle$, since zero--order H\"ormander--Mikhlin multipliers
are not bounded on $L^\infty$.  Consequently, the $L^\infty$ estimate of
Giga--Inui--Matsui \eqref{giga-28-org}, which plays a central role in the present
argument, is unavailable in the averaged setting.  In this sense, Tao's example can be interpreted as a setting in which the specific
$L^\infty$--estimate that drives our approach fails: the multipliers $m_1(D)$ and
$m_2(D)$ obstruct any direct control of the averaged product in $L^\infty$, so the
mild formulation does not yield the short--time $L^\infty$ bounds that we use to
implement a frequency--localized Picard/diagnostic construction for $\A^R u$.

\smallskip

Finally, while nonlinear effects could in principle dominate viscous dissipation and
lead to an increasingly rapid transfer of activity to small scales \cite{TaoLoc}, this
possibility is excluded on sufficiently short time intervals as long as the solution
remains bounded.  For bounded initial data, such $L^\infty$ control is available
locally in time: by a classical result of Giga--Inui--Matsui \cite{giga}, the solution
remains bounded on a time interval whose length is comparable to $\|\u_0\|_\infty^{-2}$.
This local control enables us to work with the mild formulation and to implement the
Picard iteration scheme for $\A^R\u$ under the turbulent tail hypothesis \eqref{turbulent}, thereby
contradicting any blow--up scenario consistent with \eqref{turbulent}.

\smallskip

{\em Organization of the paper.}
The paper is organized as follows. In Section~2 we prove that any Leray solution
satisfying the turbulence condition \eqref{turbulent} is bounded (and hence smooth)
on its interval of existence; the proof is based on a frequency decomposition,
a time--localized diagnostic Picard iteration for $\A^R\u$, and $L^\infty$ estimates
of Giga--Inui--Matsui type combined with Bernstein inequalities. In Section~3 we
refine this argument and show that it suffices to impose a weakened turbulence
condition only on a terminal time layer near a putative blow--up time, exploiting
high--frequency heat--flow decay on time scales comparable to $R^{-(2-\epsilon)}$, $\epsilon>0$.

\section{Turbulent solutions are bounded}

Before turning to the proof of Theorem~\ref{main}, we collect the basic analytic
estimates that will be used throughout the paper.

\subsection*{Basic inequalities}

We begin by recalling the classical existence theory of Leray together with the
energy inequality.

\begin{theorem}[Leray existence and energy inequality]\label{Leray}
For every divergence--free initial datum $\u_0\in L^2(\R^3)$, there exists a global
weak solution $\u$ to \eqref{N-S} (a Leray solution) such that
\[
\u\in L^\infty_{\mathrm{loc}}\bigl([0,\infty);L^2(\R^3)\bigr)
\cap L^2_{\mathrm{loc}}\bigl([0,\infty);H^1(\R^3)\bigr),
\qquad \u|_{t=0}=\u_0,
\]
and for almost every $t\ge 0$ the energy inequality holds:
\begin{equation}\label{energy}
\|\u(t)\|_{L^2(\R^3)}^2 + 2\int_0^t \|\nabla \u(s)\|_{L^2(\R^3)}^2\,ds
\;\le\; \|\u_0\|_{L^2(\R^3)}^2.
\end{equation}
\end{theorem}

\smallskip

Next we recall the $L^\infty$ estimate of Giga--Inui--Matsui \cite[(2.8)]{giga}.
There exists a universal constant $C_G>0$ such that for all $t>0$ and all
$f\in L^\infty(\R^3;\R^{3})$,
\begin{equation}\label{giga-28-org}
\bigl\|\nabla\!\cdot e^{t\Delta}\PP f\bigr\|_{L^\infty(\R^3)}
\le C_G\, t^{-1/2}\, \|f\|_{L^\infty(\R^3)} .
\end{equation}

We shall use the following uniform variant, valid after inserting an additional
high--frequency cutoff.

\begin{theorem}[Uniform Giga--Inui--Matsui estimate with cutoff]\label{thm:giga-cutoff}
Let $\A^R$ be defined as in \eqref{AR}. Then there exists a constant $C_D>0$,
depending only on $\phi$ (and on the dimension), such that for all $t>0$, all $R\ge 1$,
and all $f\in L^\infty(\R^3;\R^{3})$,
\begin{equation}\label{giga-28}
\bigl\|\nabla\!\cdot e^{t\Delta} \A^R \PP f\bigr\|_{L^\infty(\R^3)}
\le C_D\, t^{-1/2}\, \|f\|_{L^\infty(\R^3)}.
\end{equation}
In particular, the constant $C_D$ is independent of $R$.
\end{theorem}

\begin{proof}
We use two elementary facts:

\smallskip
\noindent
(i) $\A^R$ is a Fourier multiplier with symbol $m_R(\mxi):=1-\phi(|\mxi|/R)$, hence it
commutes with $e^{t\Delta}$, with $\nabla$, and with $\PP$;

\smallskip
\noindent
(ii) $\A^R$ is bounded on $L^\infty$ with an operator norm independent of $R$.

\smallskip
To justify (ii), note that $\A_R$ is convolution with the kernel
\[
K_R(\mx):=R^3 K(R\mx),\qquad
K:=\mathcal{F}^{-1}\!\bigl(\phi(|\mxi|)\bigr)\in\mathcal{S}(\R^3).
\]
In particular, $K\in L^1(\R^3)$ and $\|K_R\|_{L^1}=\|K\|_{L^1}$ for all $R\ge 1$.
Hence, by Young's inequality,
\[
\|\A_R g\|_{L^\infty} \le \|K_R\|_{L^1}\,\|g\|_{L^\infty}
= \|K\|_{L^1}\,\|g\|_{L^\infty},
\]
so $\|\A_R\|_{L^\infty\to L^\infty}\le \|K\|_{L^1}$ uniformly in $R$, and therefore
\[
\|\A^R g\|_{L^\infty}
\le \|g\|_{L^\infty}+\|\A_R g\|_{L^\infty}
\le (1+\|K\|_{L^1})\,\|g\|_{L^\infty}.
\]

Now apply \eqref{giga-28-org} to $\A^R f$:
\[
\bigl\|\nabla\!\cdot e^{t\Delta}\PP(\A^R f)\bigr\|_{L^\infty}
\le C_G t^{-1/2}\|\A^R f\|_{L^\infty}
\le C_G(1+\|K\|_{L^1}) t^{-1/2}\|f\|_{L^\infty}.
\]
Since all operators involved are Fourier multipliers, we have
\[
\nabla\!\cdot e^{t\Delta}\A^R\PP f
=\A^R\bigl(\nabla\!\cdot e^{t\Delta}\PP f\bigr)
=\nabla\!\cdot e^{t\Delta}\PP(\A^R f),
\]
and the desired estimate \eqref{giga-28} follows with
$C_D:=C_G(1+\|K\|_{L^1})$, which is independent of $R$.
\end{proof}

\smallskip

We also record Bernstein inequalities adapted to the low-- and high--frequency
cutoffs (see e.g.\ \cite[(26)]{TaoLoc}).

\begin{lemma}[Bernstein inequalities]\label{lem:bernstein}
Let $\A_R=\phi(|D|/R)$ and $\A^R=I-\A_R$ be as in \eqref{AR}. Then:

\begin{itemize}
\item For every $g\in L^2(\R^3)$ and every $R\ge 1$,
\begin{equation}\label{A-RA-R}
\|\A_R g\|_{L^\infty(\R^3)} \;\lesssim\; R^{3/2}\,\|g\|_{L^2(\R^3)},
\end{equation}
where the implicit constant depends only on $\phi$ (and the dimension).
In particular, for a Leray solution $\u$,
\[
\|\A_R \u(t)\|_{L^\infty}
\lesssim R^{3/2}\|\u(t)\|_{L^2}
\le R^{3/2}\|\u_0\|_{L^2},
\]
where the last inequality follows from \eqref{energy}.

\item Let $1\le p\le \infty$ and $k\in\N$ satisfy $k>\frac{3}{p}$. Then for every
$g$ with $\nabla^k g\in L^p(\R^3)$ and every $R\ge 1$,
\begin{equation}\label{A+RA+R}
\|\A^R g\|_{L^\infty(\R^3)}
\;\le\; C_{k,p}\, R^{-k+3/p}\,\|\nabla^k g\|_{L^p(\R^3)},
\end{equation}
where the constant $C_{k,p}$ depends only on $k,p,\phi$ (and the dimension), but is
independent of $R$.
\end{itemize}
\end{lemma}

\smallskip

Finally, we will use standard decay estimates for the heat semigroup with
high--frequency localization (see e.g.\ \cite[Lemma 2.2]{TaoLoc}).

\begin{lemma}[Heat flow decay with high--frequency cutoff]\label{lem:heat-decay}
Let $\u_0\in L^2(\R^3)$ and let $\A^R=I-\phi(|D|/R)$ be as above. Then for all $t>0$
and $R\ge 1$,
\begin{equation}\label{0-term}
\|e^{t\Delta}\A^R \u_0\|_{L^\infty(\R^3)}
\;\leq\; C_e t^{-3/4}e^{-ctR^2}\,\|\u_0\|_{L^2(\R^3)},
\end{equation}
for some absolute constants $c>0$ and $C_e$.  In the sequel, to simplify notation, we fix $c=1$ which can be done by adjusting the cutoff radius.

Moreover, if $\u_0\in L^\infty(\R^3)$, then
\begin{equation}\label{0-term-Linfty}
\|e^{t\Delta}\A^R \u_0\|_{L^\infty(\R^3)}
\;\le\; \|\A^R \u_0\|_{L^\infty(\R^3)}.
\end{equation}
\end{lemma}


\subsection*{Contradiction setup}

Existence of Leray solutions to \eqref{N-S} is classical since the work of
Leray~\cite{ler}, where it is shown (via retarded approximations) that one obtains a
global weak solution satisfying the energy inequality.  It remains open whether such
solutions can blow up in finite time.  In this work we show that, under the
turbulence hypothesis \eqref{turbulent} at a sufficiently large cutoff scale $R$,
this scenario cannot occur.

\smallskip

Arguing by contradiction, let $\u$ be a Leray solution corresponding to the
divergence--free Schwartz initial data $\u_0$, and define the maximal time of
boundedness by
\begin{equation}\label{blowup}
T^* := \sup\Bigl\{\,T>0:\ \u\in L^\infty\bigl((0,T)\times\R^3\bigr)\Bigr\}\in(0,\infty].
\end{equation}

Our goal is to show that \eqref{blowup} is incompatible with \eqref{turbulent}.
More precisely, assuming \eqref{turbulent} on $(0,T^*)$, we will obtain a uniform
bound on the high--frequency component $\A^R\u(t)$ for all $t<T^*$.  Together with
the low--frequency Bernstein bound \eqref{A-RA-R} (and the energy inequality
\eqref{energy}), this yields a uniform $L^\infty$ bound for $\u$ on $(0,T^*)$, which
contradicts \eqref{blowup}.

\smallskip

The key step is a time--localized Picard iteration designed to reconstruct $\A^R\u$
as a fixed point.  Throughout we interpret the expression
\[
\frac{\u(s)\otimes\u(s)}{\|\A^R\u(s)\|_\infty^2}
\]
as $0$ whenever $\|\A^R\u(s)\|_\infty=0$ (which is consistent with the subsequent
bounds).

\begin{theorem}\label{picard-convergence}
Assume that \eqref{turbulent} holds on $(0,T^*)$ for some $k\ge 2$, $p\in[1,\infty]$,
$\epsilon>0$, and some $R\ge 1$.  Suppose moreover that
\begin{equation}\label{smallness}
C_D C_{k,p}\,\sqrt{T^*}\,R^{-\epsilon} < \frac18,
\end{equation}
where $C_D$ is the constant from \eqref{giga-28} and $C_{k,p}$ is the constant from
\eqref{A+RA+R}.

Define $(\w_n)_{n\ge 0}$ by
\[
\w_0(t)\equiv 0,
\]
and for $n\ge 1$,
\begin{equation}\label{Picard}
\w_n(t)
=
e^{t\Delta}\A^R\u_0
- \int_0^t \nabla\!\cdot e^{(t-s)\Delta}\,
\A^R\PP\Big(
\frac{\u(s)\otimes\u(s)}{\|\A^R\u(s)\|_\infty^2}
\Big)\,\|\w_{n-1}(s)\|_\infty^2\,ds.
\end{equation}
Then the sequence $(\w_n)$ converges in $L^\infty\bigl((0,T^*)\times\R^3\bigr)$ to a
limit $\w\in L^\infty\bigl((0,T^*)\times\R^3\bigr)$ satisfying, for all $t\in(0,T^*)$,
\begin{equation}\label{Picard-limit}
\w(t)
=
e^{t\Delta}\A^R\u_0
- \int_0^t \nabla\!\cdot e^{(t-s)\Delta}\,
\A^R\PP\Big(
\frac{\u(s)\otimes\u(s)}{\|\A^R\u(s)\|_\infty^2}
\Big)\,\|\w(s)\|_\infty^2\,ds.
\end{equation}
Moreover, $\w(t)=\A^R\u(t)$ for every $t\in(0,T^*)$.
\end{theorem}

\begin{proof}
\emph{Step 1: uniform boundedness of the iterates.}
We claim that for all $n\ge 1$ and all $t\in(0,T^*)$,
\begin{equation}\label{wn-upper}
\|\w_n(t)\|_\infty \le 2\,\|\A^R\u_0\|_\infty .
\end{equation}
For $n=1$ this is immediate from \eqref{Picard}.  Assume \eqref{wn-upper} holds for
$\w_{n-1}$.  Using \eqref{Picard}, the cutoff estimate \eqref{giga-28}, and the
turbulence relation \eqref{turbulent}, we obtain
\begin{align*}
\|\w_n(t)\|_\infty
&\le \|e^{t\Delta}\A^R\u_0\|_\infty
 +\int_0^t \bigl\|\nabla\!\cdot e^{(t-s)\Delta}\A^R\PP F(s)\bigr\|_\infty\,
        \|\w_{n-1}(s)\|_\infty^2\,ds \\
&\le \|\A^R\u_0\|_\infty
 + C_D\int_0^t (t-s)^{-1/2}\,\|F(s)\|_\infty\,\|\w_{n-1}(s)\|_\infty^2\,ds,
\end{align*}
where
\[
F(s):=\frac{\u(s)\otimes\u(s)}{\|\A^R\u(s)\|_\infty^2}.
\]
Since $\|F(s)\|_\infty=\|\u(s)\|_\infty^2/\|\A^R\u(s)\|_\infty^2$, the turbulence
assumption gives
\[
\|F(s)\|_\infty
\le
\frac{R^{k-3/p-\epsilon}}{\|\nabla^k\u_0\|_p}.
\]
Using the inductive bound $\|\w_{n-1}(s)\|_\infty\le 2\|\A^R\u_0\|_\infty$ and
$\int_0^t (t-s)^{-1/2}\,ds = 2\sqrt{t}$, we get
\begin{align*}
\|\w_n(t)\|_\infty
&\le \|\A^R\u_0\|_\infty
 + 4C_D\,(2\sqrt{t})\,
    \frac{R^{k-3/p-\epsilon}}{\|\nabla^k\u_0\|_p}\,\|\A^R\u_0\|_\infty^2 \\
&\le \|\A^R\u_0\|_\infty
 + 8C_D\sqrt{t}\,R^{-\epsilon}\,C_{k,p}\,\|\A^R\u_0\|_\infty
\le 2\|\A^R\u_0\|_\infty,
\end{align*}
where in the last line we used Bernstein inequality \eqref{A+RA+R} applied to $g=\u_0$ and the
smallness condition \eqref{smallness}.  This proves \eqref{wn-upper}.

\smallskip
\emph{Step 2: contraction in $L^\infty((0,T^*)\times\R^3)$.}
Subtracting \eqref{Picard} at levels $n$ and $n-1$ and arguing as above yields, for
$t\in(0,T^*)$,
\begin{align*}
\|\w_n(t)-\w_{n-1}(t)\|_\infty
&\le C_D \int_0^t (t-s)^{-1/2}\,\|F(s)\|_\infty\,
   \bigl|\|\w_{n-1}(s)\|_\infty^2-\|\w_{n-2}(s)\|_\infty^2\bigr|\,ds \\
&\le C_D \int_0^t (t-s)^{-1/2}\,\|F(s)\|_\infty\,
  \bigl(\|\w_{n-1}(s)\|_\infty+\|\w_{n-2}(s)\|_\infty\bigr)
  \|\w_{n-1}(s)-\w_{n-2}(s)\|_\infty\,ds \\
&\le 4C_D\,\|\A^R\u_0\|_\infty \int_0^t (t-s)^{-1/2}\,\|F(s)\|_\infty\,
   \|\w_{n-1}(s)-\w_{n-2}(s)\|_\infty\,ds.
\end{align*}
Taking the supremum in $t\in(0,T^*)$ and using the bound on $\|F(s)\|_\infty$ gives
\[
\|\w_n-\w_{n-1}\|_{L^\infty((0,T^*)\times\R^3)}
\le
8C_D\sqrt{T^*}\,\frac{R^{k-3/p-\epsilon}}{\|\nabla^k\u_0\|_p}\,\|\A^R\u_0\|_\infty\,
\|\w_{n-1}-\w_{n-2}\|_{L^\infty((0,T^*)\times\R^3)}.
\]
By Bernstein inequality \eqref{A+RA+R} for $\A^R\u_0$ we arrive at
\[
\|\w_n-\w_{n-1}\|_{L^\infty((0,T^*)\times\R^3)}
\le
\Bigl(8C_D C_{k,p}\sqrt{T^*}\,R^{-\epsilon}\Bigr)\,
\|\w_{n-1}-\w_{n-2}\|_{L^\infty((0,T^*)\times\R^3)}.
\]
The factor in parentheses is $<1$ by \eqref{smallness}, hence $(\w_n)$ is Cauchy in
$L^\infty((0,T^*)\times\R^3)$ and converges to some
$\w\in L^\infty((0,T^*)\times\R^3)$.  Passing to the limit in \eqref{Picard} yields
\eqref{Picard-limit}.

\smallskip
\emph{Step 3: identification of the limit with $\A^R\u$.}
Since $\u$ is bounded on every compact subinterval of $(0,T^*)$, it is smooth there
and satisfies the mild formulation.  Applying $\A^R$ to the mild form yields, for
$t\in(0,T^*)$,
\begin{equation}\label{mild-high}
\A^R \u(t)
=
e^{t\Delta}\A^R\u_0
- \int_0^t \nabla\!\cdot e^{(t-s)\Delta}\,
\A^R\PP\Big(
\frac{\u(s)\otimes\u(s)}{\|\A^R\u(s)\|_\infty^2}
\Big)\,\|\A^R\u(s)\|_\infty^2\,ds.
\end{equation}
Subtracting \eqref{mild-high} from \eqref{Picard-limit} and setting
\[
E(t):=\|\w(t)-\A^R\u(t)\|_\infty,
\qquad
G(s,t):=\bigl\|\nabla\!\cdot e^{(t-s)\Delta}\A^R\PP F(s)\bigr\|_\infty,
\]
we obtain for $t\in(0,T^*)$,
\begin{align*}
E(t)
&\le \int_0^t G(s,t)\,
\bigl|\|\w(s)\|_\infty^2-\|\A^R\u(s)\|_\infty^2\bigr|\,ds \\
&\le \int_0^t G(s,t)\,
\bigl(\|\w(s)\|_\infty+\|\A^R\u(s)\|_\infty\bigr)\,E(s)\,ds.
\end{align*}
0Fix $T\in(0,T^*)$ and set
\[
B_T:=\sup_{0<s<T}\bigl(\|\w(s)\|_\infty+\|\A^R\u(s)\|_\infty\bigr)<\infty.
\]
Using \eqref{giga-28} and \eqref{turbulent} as above, we have
\[
G(s,t)\le C_D (t-s)^{-1/2}\,\|F(s)\|_\infty
\le C_D (t-s)^{-1/2}\,\frac{R^{k-3/p-\epsilon}}{\|\nabla^k\u_0\|_p}.
\]
Hence, for $t\in(0,T)$,
\[
E(t)\le C_D\,\frac{R^{k-3/p-\epsilon}}{\|\nabla^k\u_0\|_p}\,B_T
\int_0^t (t-s)^{-1/2}E(s)\,ds.
\]
Consequently,
\[
E(t)\le 2C_T\sqrt{t}\,\sup_{0<r<t}E(r),
\qquad t\in(0,T),
\]
where $C_T:= C_D\,\frac{R^{k-3/p-\epsilon}}{\|\nabla^k\u_0\|_p}\,B_T$.
Choose $\delta>0$ such that $2C_T\sqrt{\delta}<1$. Then $\sup_{0<r<\delta}E(r)=0$,
hence $E\equiv0$ on $(0,\delta)$.

Now fix any $t_0\in(0,T^*)$ and apply the same argument to the time--shifted functions
$s\mapsto \w(t_0+s)$ and $s\mapsto \A^R\u(t_0+s)$ on $s\in(0,\delta)$ (the estimates
are translation--invariant).  Iterating finitely many times covers any compact
subinterval of $(0,T^*)$, and we conclude that $E(t)\equiv 0$ on $(0,T^*)$, i.e.
\[
\w(t)=\A^R\u(t)\qquad\text{for all }t\in(0,T^*).
\]
This completes the proof.
\end{proof} Direct corollary of the previous theorem is boundedness of turbulent solution on $(0,T^*)$ contradicting \eqref{blowup}. 
\begin{corollary}[No blow--up at $T^*$ under the turbulence regime]\label{cor:no-blowup}
Assume that the turbulence condition \eqref{turbulent} holds on $(0,T^*)$ for some
parameters $k\ge 2$, $p\in[1,\infty]$, $\epsilon>0$, and some cutoff scale $R\ge 1$.
Assume moreover that the smallness condition \eqref{smallness} is satisfied.
Then $T^*$ cannot be a blow--up time; equivalently, $\u$ remains bounded on
$(0,T^*)\times\R^3$, and hence $T^*=\infty$.
\end{corollary}

\begin{proof}
If $T^*=\infty$ there is nothing to prove, so assume for contradiction that
$T^*<\infty$.

By Theorem~\ref{picard-convergence} we have the uniform bound
\[
\sup_{0<t<T^*}\|\A^R\u(t)\|_{L^\infty(\R^3)}
=\sup_{0<t<T^*}\|\w(t)\|_{L^\infty(\R^3)}
\le 2\|\A^R\u_0\|_{L^\infty(\R^3)}<\infty.
\]

On the other hand, by the low--frequency Bernstein inequality \eqref{A-RA-R} and the
energy inequality \eqref{energy}, for all $t\in(0,T^*)$,
\[
\|\A_R\u(t)\|_{L^\infty(\R^3)}
\;\lesssim\; R^{3/2}\,\|\u(t)\|_{L^2(\R^3)}
\;\le\; R^{3/2}\,\|\u_0\|_{L^2(\R^3)}.
\]
Here the implicit constant depends only on the cutoff function $\phi$ and the
dimension.

Combining the decomposition $\u=\A_R\u+\A^R\u$ with the previous bounds gives
\[
\sup_{0<t<T^*}\|\u(t)\|_{L^\infty(\R^3)}
\le
\sup_{0<t<T^*}\|\A_R\u(t)\|_{L^\infty}
+
\sup_{0<t<T^*}\|\A^R\u(t)\|_{L^\infty}
<\infty.
\]
Thus $\u\in L^\infty\bigl((0,T^*)\times\R^3\bigr)$, contradicting the definition of
$T^*$ in \eqref{blowup} as the maximal time of boundedness (under the assumption
$T^*<\infty$). Therefore $T^*$ cannot be a finite blow--up time.
\end{proof}

\section{An improvement of the turbulence condition}

In the previous section we imposed the turbulence condition \eqref{turbulent} on the
whole interval $(0,T^*)$, where $T^*$ denotes the maximal time of boundedness.
The aim of this section is twofold: we show that it is enough to assume a
turbulence--type condition only in a \emph{terminal time layer} near the putative
blow--up time, and we obtain a better smallness factor in the Picard/diagnostic
argument by exploiting the exponential damping of high frequencies over time
intervals of length $\gtrsim R^{-(2-\epsilon)}$.

\smallskip

Let $T^*$ be defined as in \eqref{blowup}, and let $0<T\le T^*$.
Fix parameters $k\ge 2$, $p\in[1,\infty]$, and $\epsilon\in(0,2)$.
Given $R\ge 1$, set the \emph{high--frequency dissipation time scale}
\begin{equation}\label{def:deltaR}
\delta_R := R^{-(2-\epsilon)}.
\end{equation}
We consider the following \emph{terminal turbulence} hypothesis: there exists
$\tau\in(0,T]$ with $\tau\ge \delta_R$ such that
\begin{equation}\label{turbulent-1}
\|\u(t)\|_\infty
\le
\frac{R^{\frac{k-3/p-\epsilon+1}{2}}}{\|\nabla^k \u(T-\tau)\|^{1/2}_{{p}}}\,
\|\A^R\u(t)\|_\infty,
\qquad t\in(T-\tau,T).
\end{equation} Compared to \eqref{turbulent}, the power of $R$ on the right--hand side is larger by
$R^{1/2}$, i.e.\ \eqref{turbulent-1} is a \emph{weaker} pointwise relation.
However, as we explain below, after restricting nonlinear interactions to the short
time scale $\delta_R$ this extra factor $R^{1/2}$ is precisely compensated by the
gain $\sqrt{\delta_R}=R^{-(1-\epsilon/2)}$, leading to an effective smallness factor
$R^{-\epsilon/2}$ in the iteration.

\smallskip

We prove that a Leray solution cannot satisfy \eqref{turbulent-1} in a neighborhood of
a blow--up time.

\begin{theorem}\label{T-main-2}
Assume that $T^*<\infty$ and fix $\epsilon\in(0,2)$, $k\ge \max\{3/p,2\}$, and $p\in[1,\infty]$.
Let $T=T^*$ and suppose that there exists $\tau\in(0,T^*]$ and $R\ge 1$ such that
$\tau\ge \delta_R$ (with $\delta_R$ as in \eqref{def:deltaR}) and the terminal
turbulence condition \eqref{turbulent-1} holds on $(T^*-\tau,T^*)$.

Let $C_D$ be the constant from \eqref{giga-28}, let $C_{k,p}$ be the constant from
\eqref{A+RA+R}, let $C_e$ be the implicit constant in \eqref{0-term} (with $c=1$),
and let $C_{\PP}:=\sup_{\mxi\neq 0}\|\PP(\mxi)\|_{\mathrm{op}}\le 1$.  Assume that $R$ is
so large that, with $\u_{in}:=\u(T^*-\tau)$,
\begin{equation}\label{ass:Rlarge}
C_e\,R^{\frac34(2-\epsilon)}e^{-R^\epsilon}\,\|\u_{in}\|_{L^2}
\;+\;\frac{C_{\PP}}{4\pi^2}\,R^{2-\epsilon}e^{-R^\epsilon/2}\,\|\u_{in}\|_{L^2}^2
\;\le\; \|\A^R\u_{in}\|_{L^\infty},
\end{equation}
and
\begin{equation}\label{ass:small}
C_D\,C_{k,p}\,R^{-\epsilon/2}<\frac18.
\end{equation}
Then \eqref{turbulent-1} cannot hold for $T=T^*$, i.e.\ a Leray solution cannot be
turbulent in the sense of \eqref{turbulent-1} in a neighborhood of the blow--up time.
\end{theorem}

\begin{proof}
Assume for contradiction that \eqref{turbulent-1} holds with $T=T^*$.
By time translation we may assume
\[
T^*-\tau=0,\qquad \u_{in}=\u(0),
\qquad\text{and}\qquad T^*=\tau,
\]
so that \eqref{turbulent-1} holds for all $t\in(0,T^*)$.
We keep the notation $\delta_R=R^{-(2-\epsilon)}$.

\smallskip
\noindent\emph{Diagnostic Picard sequence.}
Define $\w_0(t)\equiv 0$.  For $n\ge 1$ we set
\begin{equation}\label{Picard-1}
\begin{aligned}
\w_n(t)
:=\;&e^{t\Delta}\A^R\u(0)
-\int_{0}^{(t-\delta_R)_+} \nabla\!\cdot e^{(t-s)\Delta}\,\A^R\PP(\u\otimes\u)(s)\,ds\\
&\;-\int_{(t-\delta_R)_+}^{t} \nabla\!\cdot e^{(t-s)\Delta}\,
\A^R\PP\Big(\frac{\u(s)\otimes\u(s)}{\|\A^R\u(s)\|_\infty^2}\Big)\,
\|\w_{n-1}(s)\|_\infty^2\,ds,
\end{aligned}
\end{equation}
where $(t-\delta_R)_+:=\max\{t-\delta_R,0\}$, and as before we interpret the integrand
as $0$ whenever $\|\A^R\u(s)\|_\infty=0$.

The structure is the same as in Theorem~\ref{picard-convergence}, except that the
Duhamel integral is split into a \emph{far part} ($t-s\ge \delta_R$), which is
controlled using exponential high--frequency damping, and a \emph{near part}
($t-s\le \delta_R$), where we use \eqref{turbulent-1}.

\smallskip
\noindent\emph{Step 1: estimate of the far term for $\w_1$.}
For $t\le \delta_R$ the definition \eqref{Picard-1} reduces to the same form as in
Theorem~\ref{picard-convergence} and causes no difficulty.  Assume $t>\delta_R$.
Write
\[
\w_1(t)=e^{t\Delta}\A^R\u(0)\;+\;I_{\mathrm{far}}(t),
\qquad
I_{\mathrm{far}}(t)
:=-\int_{0}^{t-\delta_R} \nabla\!\cdot e^{(t-s)\Delta}\,\A^R\PP(\u\otimes\u)(s)\,ds.
\]

The high--frequency heat decay \eqref{0-term} gives
\begin{equation}\label{heat:near}
\|e^{t\Delta}\A^R\u(0)\|_\infty
\lesssim t^{-3/4}e^{-tR^2}\|\u(0)\|_{L^2}
\le C_e\,\delta_R^{-3/4}e^{-\delta_R R^2}\|\u(0)\|_{L^2}
= C_e\,R^{\frac34(2-\epsilon)}e^{-R^\epsilon}\|\u(0)\|_{L^2}.
\end{equation}

To estimate $I_{\mathrm{far}}$, we bound the $L^1\to L^\infty$ norm of the kernel.
The Fourier symbol of $\nabla e^{\tau\Delta}\A^R\PP$ is
\[
m_{\tau,R}(\mxi)=i\mxi\,e^{-\tau|\mxi|^2}\bigl(1-\phi(|\mxi|/R)\bigr)\PP(\mxi),
\]
supported in $\{|\mxi|\ge R\}$.  Since $\|\PP(\mxi)\|_{\mathrm{op}}\le C_{\PP}$, we obtain
\begin{align*}
\|\nabla e^{\tau\Delta}\A^R\PP\|_{L^1\to L^\infty}
&\le C_{\PP}(2\pi)^{-3}\int_{|\mxi|\ge R}|\mxi|e^{-\tau|\mxi|^2}\,d\mxi\\
&=C_{\PP}(2\pi)^{-3}\,4\pi\int_{R}^{\infty}r^3e^{-\tau r^2}\,dr
=\frac{C_{\PP}}{4\pi^2}\,\tau^{-2}(\tau R^2+1)e^{-\tau R^2}.
\end{align*}
Using $\tau R^2+1\le 2e^{\tau R^2/2}$ yields the uniform bound
\begin{equation}\label{kernel:L1Linf}
\|\nabla e^{\tau\Delta}\A^R\PP\|_{L^1\to L^\infty}
\le \frac{C_{\PP}}{4\pi^2}\,\tau^{-2}e^{-\tau R^2/2},
\qquad \tau>0,\ R\ge 1.
\end{equation}

By the energy inequality \eqref{energy}, for a.e.\ $s$ one has
$\|\u(s)\otimes\u(s)\|_{L^1}\le \|\u(s)\|_{L^2}^2\le \|\u(0)\|_{L^2}^2$.
Therefore,
\begin{align*}
\|I_{\mathrm{far}}(t)\|_{\infty}
&\le \int_0^{t-\delta_R}
\|\nabla e^{(t-s)\Delta}\A^R\PP\|_{L^1\to L^\infty}\,
\|\u(s)\otimes\u(s)\|_{L^1}\,ds \\
&\le \frac{C_{\PP}}{4\pi^2}\,\|\u(0)\|_{L^2}^2
\int_{\delta_R}^{\infty}\tau^{-2}e^{-\tau R^2/2}\,d\tau.
\end{align*}
Changing variables $u=\tfrac12 R^2\tau$ gives
\[
\int_{\delta_R}^{\infty}\tau^{-2}e^{-\tau R^2/2}\,d\tau
=\frac{R^2}{2}\int_{R^\epsilon/2}^{\infty}u^{-2}e^{-u}\,du
\le \frac{R^2}{2}\,e^{-R^\epsilon/2}\int_{R^\epsilon/2}^{\infty}u^{-2}\,du
=R^{2-\epsilon}e^{-R^\epsilon/2},
\]
and hence
\begin{equation}\label{Ifar:bound}
\|I_{\mathrm{far}}(t)\|_{\infty}
\le \frac{C_{\PP}}{4\pi^2}\,R^{2-\epsilon}e^{-R^\epsilon/2}\,\|\u(0)\|_{L^2}^2.
\end{equation}
Combining \eqref{heat:near}--\eqref{Ifar:bound} with \eqref{ass:Rlarge} yields
\begin{equation}\label{w1:bound}
\|\w_1(t)\|_\infty \le \|\A^R\u(0)\|_{L^\infty},\qquad t\in(0,T^*).
\end{equation}

\smallskip
\noindent\emph{Step 2: uniform bounds on $(\w_n)$ and contraction.}
We claim that for all $n\ge 1$ and all $t\in(0,T^*)$,
\begin{equation}\label{wn:bound-local}
\|\w_n(t)\|_\infty \le 2\,\|\A^R\u(0)\|_{L^\infty}.
\end{equation}
For $n=1$ this is \eqref{w1:bound}.  Assume \eqref{wn:bound-local} holds for $n-1$.
Using \eqref{Picard-1}, the estimate \eqref{giga-28} on the near part, and
\eqref{turbulent-1}, we obtain for $t\in(0,T^*)$,
\begin{align*}
\|\w_n(t)\|_\infty
&\le \|\w_1(t)\|_\infty
+ C_D\int_{(t-\delta_R)_+}^{t}(t-s)^{-1/2}\,
\frac{\|\u(s)\|_\infty^2}{\|\A^R\u(s)\|_\infty^2}\,
\|\w_{n-1}(s)\|_\infty^2\,ds\\
&\le \|\w_1(t)\|_\infty
+ C_D\int_{(t-\delta_R)_+}^{t}(t-s)^{-1/2}\,
\frac{R^{k-3/p-\epsilon+1}}{\|\nabla^k\u(0)\|_{p}}\,
\|\w_{n-1}(s)\|_\infty^2\,ds.
\end{align*}
Since $\int_{(t-\delta_R)_+}^{t}(t-s)^{-1/2}\,ds\le 2\sqrt{\delta_R}
=2R^{-(1-\epsilon/2)}$, and $\|\w_{n-1}\|_\infty\le 2\|\A^R\u(0)\|_\infty$,
\[
\|\w_n(t)\|_\infty
\le \|\w_1(t)\|_\infty
+8C_D\,R^{-(1-\epsilon/2)}\,\frac{R^{k-3/p-\epsilon+1}}{\|\nabla^k\u(0)\|_{p}}\,
\|\A^R\u(0)\|_\infty^2.
\]
Using Bernstein inequality \eqref{A+RA+R} at time $0$ (note that $\u(0)$ is smooth), we have
$\|\A^R\u(0)\|_\infty\le C_{k,p}R^{-k+3/p}\|\nabla^k\u(0)\|_{p}$, hence
\[
\|\w_n(t)\|_\infty
\le \|\w_1(t)\|_\infty
+8C_D C_{k,p}\,R^{-\epsilon/2}\,\|\A^R\u(0)\|_\infty
\le 2\|\A^R\u(0)\|_\infty
\]
by \eqref{w1:bound} and \eqref{ass:small}.  This proves \eqref{wn:bound-local}.

The contraction estimate is identical: subtracting \eqref{Picard-1} at levels $n$
and $n-1$, only the near integral contributes, and using the same bounds yields
\[
\|\w_n-\w_{n-1}\|_{L^\infty((0,T^*)\times\R^3)}
\le \bigl(8C_D C_{k,p}R^{-\epsilon/2}\bigr)
\|\w_{n-1}-\w_{n-2}\|_{L^\infty((0,T^*)\times\R^3)}.
\]
By \eqref{ass:small} the factor is $<1$, hence $(\w_n)$ converges in
$L^\infty((0,T^*)\times\R^3)$ to a fixed point $\w$.

\smallskip
\noindent\emph{Step 3: identification $\w=\A^R\u$.}
As in the proof of Theorem~\ref{picard-convergence}, $\u$ is smooth on compact
subintervals of $(0,T^*)$ and satisfies the mild formulation there.
Applying $\A^R$ yields an identity for $\A^R\u$ analogous to \eqref{Picard-1} with
$\|\w_{n-1}(s)\|_\infty^2$ replaced by $\|\A^R\u(s)\|_\infty^2$.
Subtracting the two fixed point relations and arguing on short time steps
(of length $\delta_R$) shows that $\|\w(t)-\A^R\u(t)\|_\infty\equiv 0$ on $(0,T^*)$.
Consequently, $\A^R\u$ is bounded on $(0,T^*)$.

Finally, Bernstein inequality for the low frequencies \eqref{A-RA-R} together with the energy
inequality \eqref{energy} yields $\|\A_R\u(t)\|_\infty\lesssim R^{3/2}\|\u(0)\|_2$.
Therefore $\u=\A_R\u+\A^R\u$ is bounded on $(0,T^*)\times\R^3$, contradicting the
definition of $T^*$ and \eqref{blowup}.  This completes the proof.
\end{proof}

\section*{Declarations}

\subsection*{Conflict of interest}
The author declares that there are no conflicts of interest.

\subsection*{Data availability}
Data sharing is not applicable to this article as no datasets were generated or
analysed during the current study.

\subsection*{Funding}
This work was supported in part by the Austrian Science Fund (FWF), project
P35508 (DOI: 10.55776/P35508).

\end{document}